\documentclass[11pt]{article}
\def\mytitle{An efficient  second-order cone programming approach 
for optimal selection in tree breeding}

\usepackage{amssymb}
\usepackage{amsthm}
\usepackage{color}
\usepackage{enumerate}
\usepackage{url}
\usepackage{graphicx}
 
\def\@themcountersep{}
\def\0{\mbox{\bf 0}}
\def\1{\mbox{\bf 1}}
\def\2{\mbox{\bf 2}}
\def\3{\mbox{\bf 3}}
\def\4{\mbox{\bf 4}}
\def\5{\mbox{\bf 5}}
\def\6{\mbox{\bf 6}}
\def\7{\mbox{\bf 7}}
\def\8{\mbox{\bf 8}}
\def\9{\mbox{\bf 9}}

\def\c{\mbox{\boldmath $c$}}

\newdimen\zhige \zhige=0pt
\def\chige#1{{\setbox\zhige\hbox{#1}\ifdim\ht\zhige=1ex\accent24 #1%
  \else\ooalign{\unhbox\zhige\crcr\hidewidth\char24\hidewidth}\fi}}

\def\e{\mbox{\boldmath $e$}}
\def\f{\mbox{\boldmath $f$}}
\def\g{\mbox{\boldmath $g$}}
\def\h{\mbox{\boldmath $h$}}

\def\l{\mbox{\boldmath $l$}}

\def\u{\mbox{\boldmath $u$}}
\def\v{\mbox{\boldmath $v$}}

\def\x{\mbox{\boldmath $x$}}
\def\y{\mbox{\boldmath $y$}}
\def\z{\mbox{\boldmath $z$}}
\def\A{\mbox{\boldmath $A$}}
\def\B{\mbox{\boldmath $B$}}

\def\F{\mbox{\boldmath $F$}}

\def\I{\mbox{\boldmath $I$}}

\def\O{\mbox{\boldmath $O$}}

\def\U{\mbox{\boldmath $U$}}

\def\KC{\mbox{$\cal K$}}

\def\PC{\mbox{$\cal P$}}

\def\Real{\mbox{$\mathbb{R}$}}

\def\SMAT{\mbox{$\mathbb{S}$}}

\def\ringaccent#1{{\accent23 #1}}

\begin{document}

\noindent \textcolor{blue}{\leaders\hrule width 0pt height 0.08cm \hfill}
\vspace{0.5cm} 

{\Large \bf \noindent \mytitle}

\vspace{0.2cm}
 \noindent
Makoto Yamashita
\footnote{
Department of Mathematical and Computing Sciences,
 Tokyo Institute of Technology, 2-12-1-W8-29 Ookayama, Meguro-ku, Tokyo
 152-8552, Japan (Makoto.Yamashita@is.titech.ac.jp).
 },
Tim J. Mullin
\footnote{
The Swedish Forestry Research Institute (Skogforsk), Box 3, S\"{a}var 918 21, Sweden; 
and 224 rue du Grand-Royal Est, QC, J2M 1R5, Canada.
}
and
Sena Safarina
\footnote{
Department of Mathematical and Computing Sciences,
 Tokyo Institute of Technology, 2-12-1-W8-29 Ookayama, Meguro-ku, Tokyo
 152-8552, Japan.
 }
\\
Submitted: June 15, 2015.

\vspace{0.3cm}

\noindent {\bf Abstract:} 


An important problem in tree breeding is optimal selection from candidate pedigree members
to produce the highest performance in seed orchards, while conserving essential genetic diversity.
The most  beneficial  members should contribute as much as possible, but such selection of orchard parents would 
reduce performance of the orchard progeny due to serious inbreeding.
To avoid inbreeding, we should include a constraint on the numerator relationship matrix to keep a group coancestry under an appropriate threshold. 
Though an SDP (semidefinite programming) approach proposed by Pong-Wong and Woolliams 
gave an accurate optimal value, it required rather long computation time.

In this paper, we propose an SOCP (second-order cone programming) approach to reduce this computation time.
We demonstrate that the same solution is attained by the SOCP formulation,
but requires much less time.
Since a simple SOCP formulation is not much more efficient compared to the SDP approach,
we exploit a sparsity structure of the numerator relationship matrix, 
and formulate the SOCP constraint using Henderson's algorithm.
Numerical results show that 
the proposed SOCP approach reduced  computation time in a case study from 39,200 seconds under the SDP approach 
to less than 2 seconds.

\vspace{0.3cm} 

\noindent {\bf Keywords:} Semidefinite programming, Second-order cone programming, Tree breeding,
Optimal selection, Group coancestry, Relatedness, Genetic gain

\noindent {\bf AMS classification:} 90C22 Semidefinite programming, 90C25 Convex programming, 
92-08  Biology and other natural sciences (Computational methods).

\noindent \textcolor{blue}{\leaders\hrule width 0pt height 0.08cm \hfill}
\vspace{0.5cm}

\section{Introduction}\label{sec:introduction}

The usage of mathematical optimization approaches for tree breeders have been
increasing~\cite{AHLINDER14, MEUWISSEN97, PONG-WONG07, schierenbeck2011controlling},
since one of their purposes is to derive better performance 
from seed orchards.
When tree breeders make a plan for new seed orchards, they determine 
the contributions of candidate pedigree members so that the resultant orchard 
maximizes response to the selection.
From the viewpoint of mathematical optimization,
the simplest form of the optimal selection problem is:
\begin{eqnarray*}
\begin{array}{lcl}
\max &:&  \g^T\x \label{eq:simple} \\
\mbox{subject to} &:& \e^T \x = 1 \\ 
& & \l \le \x \le \u.
\end{array} 
\end{eqnarray*}
Here, we assume that the number of the candidate members is $m$.
The variable vector is $\x \in \Real^m$, and corresponds to the contributions of 
the candidates.
The first constraint $\e^T \x = 1$ indicates that the total contribution of candidate members is unity.
We use the vector $\e \in \Real^m$ to denote the vector of all ones,
and the superscript $T$ to denote the transpose of a vector or a matrix.
In the second constraints, $\l \in \Real^m$ and $\u \in \Real^m$ 
are element-wise lower and upper bounds on the contributions, respectively.
The performance  measure appears  in the objective function
as its coefficient $\g = (g_1, \ldots, g_m)^T$,
and the  
estimated breeding value (EBV)~\cite{LYNCH98} is often employed for $\g$.
The values $g_1, \ldots, g_m$ are calculated separately before we solve the optimal selection problem, 
so we can consider $g_1, \ldots, g_m$ as constant values.

The above problem 
is a simple linear programming problem,
therefore, a greedy method is enough to solve it. Such solution includes the candidates corresponding to
the highest EBV as much as possible, and 
it is most efficient in situations where all the pedigrees are independent.
Even if the pedigrees are independent, diversity is still an issue.
Lindgren \textit{et al.}~\cite{lindgren1989deployment} discussed a linear deployment
in which the contributions of the candidate members are proportional to their EBVs.
In practical situations, however, we cannot overlook the effects due to the relatedness that accumulates over cycles of breeding the pedigrees.

To reflect the effect of the relatedness, 
Cockerham~\cite{COCKERHAM67} 
extended the definitions of coancestry coefficients in order to include coancestry of a group.
The group coancestry on the contributions $\x$
is calculated with the formula $\frac{\x^T \A \x}{2}$, where $\A \in \Real^{m \times m}$
is the numerator relationship matrix of Wright~\cite{WRIGHT22}
(we will review a formula for the numerator relationship matrix in
 Section~\ref{sec:SDP}).
Introducing a constraint to keep group coancestry under an appropriate level $\theta \in \Real$,
Meuwissen~\cite{MEUWISSEN97} proposed a formulation of optimal contributions:
\begin{eqnarray}
\begin{array}{lcl}
\max &:& \g^T \x \label{eq:OPSEL} \\
\mbox{subject to} &:& \e^T \x = 1 \\ 
& & \frac{\x^T \A \x}{2} \le \theta \\
& & \l \le \x \le \u.
\end{array} 
\end{eqnarray}
Meuwissen developed an iterative method based on Lagrangian multipliers to solve this optimization problem,
and his method has been used in breedings, for example,  
\cite{GRUNDY98, HINRICHS11, WOOLIAMS07}.
It is a characteristic of this method that some variables $x_i$ may be fixed  to its lower or upper bounds ($l_i$ or $u_i$) during
the iterations, and Pong-Wong and Woolliams~\cite{PONG-WONG07} 
demonstrated that the method did not always 
obtain the optimal solution.

Pong-Wong and Woolliams utilized the structure of the numerator relationship matrix $\A$ to
formulate the problem~(\ref{eq:OPSEL}) into a semidefinite programming (SDP) problem.
SDP is a convex optimization problem that maximizes a linear objective function
over the constraints described as linear matrix inequalities.
The research in 1990s, for example 
\cite{HELMBERG96, KOJIMA94},
extended interior-point methods from linear programming problems to 
SDP problems. 
Based on primal-dual interior-point methods, 
software packages (like SDPA~\cite{SDP-HANDBOOK}, 
SDPARA~\cite{YAMASHITA12}, SDPA-C~\cite{yamashita2015fast},  
SDPT3~\cite{TODD99}, and SeDuMi~\cite{STURM99}) have been developed for solving SDPs.
Using the SDP formulation, Pong-Wong and Woolliams obtained the exact optimal value of~(\ref{eq:OPSEL}).
The number of candidate members discussed in~\cite{PONG-WONG07}, however, 
was limited to only small sizes,  $m \le 10$.
Ahlinder~{\it et al.}~\cite{AHLINDER14} implemented the SDP approach into a software package
called OPSEL~\cite{mullin2014opsel}\footnote{\url{http://www.skogforsk.se/opsel}} with the help of the latest version of SDPA (a high-performance SDP solver)~\cite{SDP-HANDBOOK}.
They solved large problems 
($m \ge 10,000$) that were generated from real Scots pine pedigrees and performance data, and
they also focused on flexibility and re-optimization of the SDP approach.

A main obstacle in the numerical tests of \cite{AHLINDER14}
was that the SDP approach took rather long computation time.
They reported for their case study that 
it required five-hours of computation time for a problem of the size $m = 12,000$.
Even though the SDP guarantee  the optimal solution,
the computation time is rather long for operational application and requiring significant 
truncation of the candidate list prior to optimizing the selection.

In this paper, we propose a second-order cone programming (SOCP) approach.
SOCP is a convex optimization that maximizes a linear objective function over second-order cone constraints.
SOCP can be considered as a special case of SDP,
and can be efficiently solved with interior-point methods in a similar way to SDP~\cite{SCHMIETA01, TSUCHIYA99}.
Lobo~{\it et al.}~\cite{LOBO98} discussed wide-range applications of SOCP, for example, filter design and
 truss design, and Sasakawa and Tsuchiya applied SOCP to magnetic shield design~\cite{sasakawa2003optimal}.
Alizadeh and Goldfarb~\cite{ALIZADEH03} surveyed theoretical and algorithmic
aspects of SOCP, and 
the software packages SDPT3~\cite{TODD99} and SeDuMi~\cite{STURM99} can 
solve not only SDP but also SOCP using the primal-dual interior-point methods.
In addition, ECOS~\cite{DOMAHIDI13} was also implemented recently to solve SOCP problems.

We first discuss that the proposed SOCP formulation also attains the optimal solution of~(\ref{eq:OPSEL}).
Since SOCP is a special case of SDP, we could expect that a simple SOCP formulation
would be enough to reduce the computation time. 
However, preliminary numerical tests showed that the simple formulation did not perform well.
We therefore utilize a sparsity embedded in the numerator relationship matrix
and establish a more efficient SOCP formulation.
Furthermore, we integrate Henderson's algorithm~\cite{HENDERSON76} into this formulation.
Numerical tests with the data including Scots pine showed
that the SOCP formulation with Henderson's algorithm
reduced a great amount of computation time.
For the case of $m=10,100$, we attained  a speedup of 
20,000-times compared to the SDP approach.

The rest of this paper is organized as follows.
Section~\ref{sec:SDP} describes the SDP approach of Pong-Wong and Woolliams
and discusses a simple SOCP formulation.
In Section~\ref{sec:SOCP}, we propose SOCP formulations and derive an efficient method to solve 
problem~(\ref{eq:OPSEL}) .
Section~\ref{sec:test} shows the numerical results to verify the computation time reduction 
for problems of various sizes.
Finally, Section~\ref{sec:conclusion} gives conclusions and discusses future directions.

\section{SDP formulation and simple SOCP formulation}\label{sec:SDP}

Since a principal characteristic of our problem~(\ref{eq:OPSEL}) 
is determined 
by the numerator relationship matrix $\A$, we first review a formula to evaluate 
its elements.
We then describe the SDP approach of Pong-Wong and Woolliams~\cite{PONG-WONG07},
and compare the performance of the SDP approach and a simple SOCP formulation.

\begin{figure}
\begin{minipage}{.5\textwidth}
{\small
\begin{center}
\begin{tabular}{c|c}
\hline
pedigree id & parents \\
\hline
 1 & unknown and unknown \\
 2 & unknown and unknown \\
 3 & 1 and 2 \\
 4 & 1 and 2 \\
 5 & 2 and unknown \\
 6 & 3 and 4 \\
 7 & 1 and 5 \\
 8 & 6 and 7 \\
 9 & 5 and 7 \\
\hline
\end{tabular}
\end{center}
}
\end{minipage}
\hfill
\begin{minipage}{.5\textwidth}
\begin{center}
\includegraphics[scale=0.5]{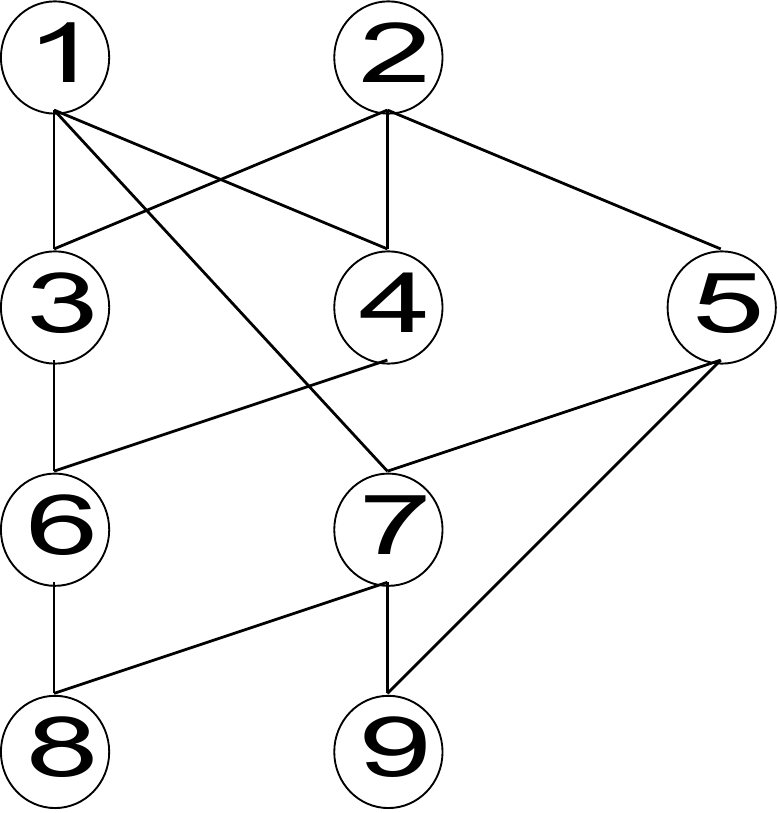}
\end{center}
\end{minipage}
\caption{An example of pedigree heredity and its diagram.}\label{fig:tree}
\end{figure}

To evaluate the elements of  the numerator relationship matrix $\A$,
we separate the set of pedigree candidate members 
$\PC := \{1, 2, \ldots, m  \}$ into the three disjoint groups:
\begin{eqnarray*}
\PC = \PC_0 \cup \PC_1 \cup \PC_2,
\end{eqnarray*}
where
\begin{eqnarray*}
\left\{
\begin{array}{rcl}
\PC_0 &=& \{ i \in \PC : \mbox{both parents} \ p(i) \ \mbox{and} \ q(i) 
\ \mbox{are unknown} \} \\
\PC_1 &=& \{ i \in \PC : \mbox{one parent} \ p(i)  
\ \mbox{is known and the other parent} \ q(i) \ \mbox{is unknown} \} \\
\PC_2 &=& \{ i \in \PC : \mbox{both parents} \ p(i) \ \mbox{and} \ q(i) 
\ \mbox{are known} \} .
\end{array}
\right.
\end{eqnarray*}
Figure~\ref{fig:tree} gives an example of pedigree  with $m=9$ members
and illustrates its genealogical chart.
In this example, $\PC_0 = \{1, 2\}, \PC_1 = \{5\}, \PC_2 = \{3, 4, 6, 7, 8, 9\}$,
and the parents of the 8th member are $p(8)  = 7$ and $q(8) = 6$.
We use a convention $p(i) = 0$ or $q(i) = 0$ if the parent $p(i)$ or $q(i)$ is unknown, respectively,
and we can assume 
$i > p(i) \ge q(i)$  for all $i \in \PC$ without loss of generality.

The numerator relationship matrix $\A$ was defined by Wright~\cite{WRIGHT22},
and its simplified formula was devised in \cite{HENDERSON76}. 
The formula of~\cite{HENDERSON76} gives the elements $A_{11}, \ldots, A_{nn}$ in a recursive style:
\begin{eqnarray*}
\left\{\begin{array}{ll}
A_{ij} = A_{ji} = \frac{A_{j,p(i)} + A_{j,q(i)}}{2}  & \mbox{for} \
i = 1, \ldots, m, \ j = 1,\ldots,i-1 \\
A_{ii} = 1 + \frac{A_{p(i),q(i)}}{2}  & \mbox{for} \
i = 1, \ldots, m,
\end{array}\right.
\end{eqnarray*}
where we use a convention $A_{pq} = 0$ if $p = 0$ or $q = 0$.
When we apply this calculation to the example of Figure~\ref{fig:tree}, 
we obtain the corresponding matrix $\A$ as follow:
{\small
\begin{eqnarray}
\A = 
\frac{1}{32}
\left(\begin{array}{rrrrrrrrr}
32  &   0  &  16  &  16  &    0  &   16  &   16  &   16  &    8 \\
    0  &   32 &    16  &   16  &   16  &   16  &    8 &    12  &   12 \\
    16   &  16  &   32  &   16   &   8  &   24  &   12   &  18  &   10 \\
    16  &   16 &    16 &    32   &   8  &   24 &    12 &    18 &    10 \\
     0  &   16 &     8  &    8  &   32  &    8  &   16  &   12  &   24 \\
    16  &   16  &   24  &   24  &    8  &   40  &   12 &    26 &    10\\
    16  &    8  &   12  &   12   &  16  &   12  &   32 &    22 &    24 \\
    16  &   12  &   18  &   18  &   12  &   26  &   22  &   38 &    17 \\
     8  &   12 &    10 &    10  &   24 &    10 &    24 &    17 &    40
\end{array}\right).
\label{eq:A}
\end{eqnarray}
} 

Pong-Wong and Woolliams~\cite{PONG-WONG07} formulated the 
problem~(\ref{eq:OPSEL}) into an SDP problem.
A standard form of SDP  is given by
\begin{eqnarray}
\begin{array}{lclcrcc}
\max &:& & & \sum_{k=1}^m c_k z_k \label{eq:SDP} \\
\mbox{subject to} &:& \F_0 &-& \sum_{k=1}^m \F_k z_k &\in &\SMAT_+^n. 
\end{array} 
\end{eqnarray}
We use $\SMAT^n$ to denote the space of $n \times n$ symmetric matrices,
and $\SMAT_+^n \subset \SMAT^n $ to denote the space of positive semidefinite 
matrices of dimension $n$. The variables are $z_1, \ldots, z_m \in \Real$, 
and the input data
are $c_1, \ldots, c_m \in \Real$ and $\F_0, \F_1, \ldots, \F_m \in \SMAT^n.$

The key step of Pong-Wong and Woolliams~\cite{PONG-WONG07} was the usage of the Schur complement
of a matrix block.
They noticed that the numerator relationship matrix $\A$ is always positive definite,
and they utilized this property to convert the constraint on
the group coancestry to a positive semidefinite condition on the symmetric matrix:
\begin{eqnarray}
\frac{\x^T \A \x}{2} \le \theta \quad \Leftrightarrow  \quad
\left(\begin{array}{cc}
2 \theta & \x^T \\
\x & \A^{-1}
\end{array}\right) \in \SMAT_+^{1+m}. \label{eq:Schur}
\end{eqnarray}
With this positive semidefinite constraint, they converted (\ref{eq:OPSEL}) 
into the standard SDP form.
The input vector $\c$ and the input matrices $\F_0, \F_1, \ldots, \F_m$ 
in~(\ref{eq:SDP})
are given as follow~\cite{AHLINDER14, PONG-WONG07}:
{\small
\begin{eqnarray*}
\c &=& \g, \\
\F_0 &=& \left(\begin{array}{ccccc}
1 &      &   &   &   \\
   & -1 &   &    &   \\
   &     & \mbox{Diag}(\u) & & \\
   &     &            & -\mbox{Diag}(\l) & \\
   &     &            & & \left(\begin{array}{cc} 2 \theta & \0^T \\ \0 & \A^{-1} \end{array}\right)
\end{array}\right), \\
\F_k &=& \left(\begin{array}{ccccc}
1 &      &   &   &   \\
   & -1 &   &    &   \\
   &     & \mbox{Diag}(\e_i) & & \\
   &     &            & -\mbox{Diag}(\e_i) & \\
   &     &            & & \left(\begin{array}{cc} 0 & -\e_i^T \\ \e_i & \O \end{array}\right)
\end{array}\right) \qquad \mbox{for} \quad k = 1, \ldots, m
\end{eqnarray*}
} 
We use $\e_i$ to denote the vector of all zeros except 1 at the $i$th element, 
and $\mbox{Diag}(\u)$ to denote the diagonal matrix whose diagonal elements are $\u$.
Note that the dimension of $\F_0, \F_1, \ldots, \F_m$ is $3m + 3$.
The software package OPSEL~\cite{mullin2014opsel} automatically formulates
the optimal selection into an SDP problem of this form.

Table~\ref{table:timeSDP} shows the computed optimal values and the computation time of 
Meuwissen's implementation (GENCONT)~\cite{MEUWISSEN02} and the SDP formulation. 
We executed the numerical tests using Matlab R2015a on Windows 8.1  PC
with Xeon CPU E3-1231 (3.40 GHz, 4 cores) and 8 GB memory space.
We used Windows,
since GENCONT can run only on Windows.
To solve the SDP~(\ref{eq:SDP}), we employed SDPA~\cite{SDP-HANDBOOK}.
\begin{table}
\caption{Computation time on GENCONT and SDP formulation (time in seconds)}
\label{table:timeSDP}
\begin{center}
\begin{tabular}{rl|r|r}
\hline
\multicolumn{2}{r|}{$m$ (the number of pedigree)} & 2,045 & 10,100 \\
\hline
GENCONT & optimal value & 438.56 & OOM${}^*$ \\
& time & 67.43 &  \\
\hline
SDP formulation & optimal value  & 439.12  & 47.76 \\
 (with 4 cores)& time & 70.21 & 39200.78 \\
\hline
\footnotesize{*OOM - ``out of memory"}
\end{tabular}
\end{center}
\end{table}

We observe from Table~\ref{table:timeSDP} that the SDP formulation
attained a better optimal value than GENCONT. Actually, 
as shown in~\cite{PONG-WONG07}, the optimal value of the SDP formulation
was the exact optimal value, while the Lagrangian multiplier method implemented in GENCONT 
could not guarantee the optimality.
For the large problem ($m = 10,100$), GENCONT gave up the computation,
but the SDP formulation again obtained the optimal solution.
On the other hand, the disadvantage of the SDP formulation is its computation time.
Even using the parallel computing implemented in SDPA,
the SDP formulation was slower than GENCONT for $m = 2,045$.
Furthermore, the computation time for the large problem exceeded 10 hours
even with four cores.
When we used only one core, the SDP formulation would require longer computation time than 24 hours.

To reduce the heavy computation time of the SDP formulation, we review the 
group coancestry constraint $\frac{\x^T \A \x}{2} \le \theta$.
With the positive definiteness of $\A$, we focus on the property that 
this constraint can also be described as a second-order condition.
The vector $\v \in \Real^{n_q} $ is said to satisfy the second-order condition
if $\v \in \KC^{n_q}$. 
The symbol $\KC^{n_q}$ denotes the second-order cone of dimension $n_q$:
\begin{eqnarray*}
\KC^{n_q} := 
\left\{ \v \in \Real^{n_q} : v_1 \ge \sqrt{\sum_{k=2}^{n_q} v_k^2} \right\}.
\end{eqnarray*}
The second-order cone is a special case of positive semidefinite constraint.
In fact, using the Schur complement, we can verify that
\begin{eqnarray*}
\v \in \KC^{n_q}
\quad  \Leftrightarrow \quad
\left(\begin{array}{ccccccccc}
v_1 & v_2 & v_3 & \cdots & v_n \\
v_2 & v_1 &  0   & \cdots & 0 \\
v_3 &  0 & v_1     & \cdots & 0 \\
\vdots & 0 & \vdots & \ddots  & 0 \\
v_n & 0  & 0     & \cdots & v_1 \\
\end{array}\right)
\in \SMAT_+^{n_q}.
\end{eqnarray*}

Furthermore,
since the numerator relationship matrix $\A$ is positive definite, we can apply 
the Cholesky factorization to $\A$ to obtain
the upper triangular matrix $\U$ such that  $\A = \U^T \U$.
This factorization derives a second-order cone condition that corresponds to 
the group coancestry constraint:
\begin{eqnarray}
\frac{\x^T \A \x}{2} \le \theta
\quad \Leftrightarrow \quad
\x ^T \U^T \U \x \le 2 \theta
\quad \Leftrightarrow \quad
|| \U \x || \le \sqrt{2 \theta}
\quad \Leftrightarrow \quad
\left(\begin{array}{c}
\sqrt{2 \theta} \\ \U \x
\end{array}\right) \in \KC^{1 + m}. \label{eq:SOCPconvert}
\end{eqnarray}

Since the most difficult constraint in~(\ref{eq:OPSEL}) can be expressed using
a second-order cone, it is natural to consider its SOCP formulation.
A standard form of SOCP problem is described by
\begin{eqnarray}
\begin{array}{lclcrcc}
\max &:& & & \c^T \z \label{eq:SOCP} \\
\mbox{subject to} &:& \f_0 &-& \F \z &\in& \Real_+^{n_l} \times \KC^{n_q}. 
\end{array} 
\end{eqnarray}
We use $\Real_+^{n_l} \times \KC^{n_q}$ to denote the Cartesian product
of the non-negative orthant of dimension $n_l$, 
$\Real_+^{n_l} := \{ \v \in \Real^{n_l} : v_i \ge 0 \ \mbox{for} \ i=1,\ldots,n_l\}$, 
 and the second-order cone of dimension $n_q$.
In this problem,  the variable vector is $\z \in \Real^{m}$,
and the input data are $\c \in \Real^m, \f_0 \in \Real_+^{n_l + n_q}$ 
and $\F \in \Real^{(n_l + n_q) \times m}$.
A more general SOCP form than (\ref{eq:SOCP}) can be defined so that it can handle the Cartesian product of multiple second-order cones,
but one cone is enough for the discussion in this paper.

If we formulate the optimal selection problem~(\ref{eq:OPSEL}) in a simple style, 
we obtain an SOCP problem:
\begin{eqnarray}
\begin{array}{lclcrcc}
\max &:& & &  \g^T \x \label{eq:straightSOCP} \\
\mbox{subject to} &:& 
\left(\begin{array}{r}
1 \\ -1 \\ \u \\ -\l \\ 
\hline
\sqrt{2 \theta} \\ \0
\end{array} \right)
&-&
\left(\begin{array}{r}
\e^T \\ -\e^T \\ \I \\ -\I \\ 
\hline
0 \\ \U
\end{array} \right) \x & \in &
\Real_+^{2 + 2m} \times \KC^{1+m},
\end{array}
\end{eqnarray}
We use $\I$ to denote the identity matrix of dimension $m$.
We call this formulation {\it a simple SOCP formulation}.
We emphasize that this simple SOCP formulation also gives the exact optimal value
of the optimal selection problem~(\ref{eq:OPSEL})
due to the equivalence from (\ref{eq:Schur}) and (\ref{eq:SOCPconvert});
\begin{eqnarray*}
\left(\begin{array}{cc} 2\theta & \x^T \\ \x & \A^{-1} \end{array} \right)
\in \SMAT_+^{1+m}
\quad
\Leftrightarrow
\quad
\left(\begin{array}{c} \sqrt{2\theta}  \\ \U\x  \end{array} \right)
\in \KC^{1+m}.
\end{eqnarray*}

Since SOCP is a special case of SDP, we expected that 
we could solve the SOCP formulation~(\ref{eq:straightSOCP})
faster than the SDP formulation~(\ref{eq:SDP}).
In Table~\ref{table:straightSOCP}, we compare the computation times 
of the SDP and SOCP formulations. We used ECOS~\cite{DOMAHIDI13} as the SOCP solver
for~(\ref{eq:straightSOCP}).
\begin{table}
\caption{Computation time on SDP and simple SOCP formulations (time in seconds)}
\label{table:straightSOCP}
\begin{center}
\begin{tabular}{r|r|r}
\hline
$m$ (the number of pedigree) & 2,045 & 10,100 \\
\hline
SDP formulation  (\ref{eq:SDP}) & 70.21 & 39200.78 \\
simple SOCP formulation (\ref{eq:straightSOCP}) & 0.28 & 5604.25 \\
\hline
\end{tabular}
\end{center}
\end{table}
For the small problem ($m = 2,045$), the SOCP formulation successfully reduced
the computation time from 70.21 seconds to 0.28 seconds, so its speed-up was 250-times.
However, for the large problem ($m = 10,100$), the speed-up was limited to 7-times.
When we consider the computational complexity, the simple SOCP formulation would be slower 
than the SDP formulation for further large problems.
On contrary to the fact that SOCP is a special case of SDP, 
the simple SOCP formulation did not seem promising.

\section{Efficient formulation based on SOCP}\label{sec:SOCP}

\begin{figure}
\begin{minipage}{.5\textwidth}
\begin{center}
\includegraphics[width=200pt]{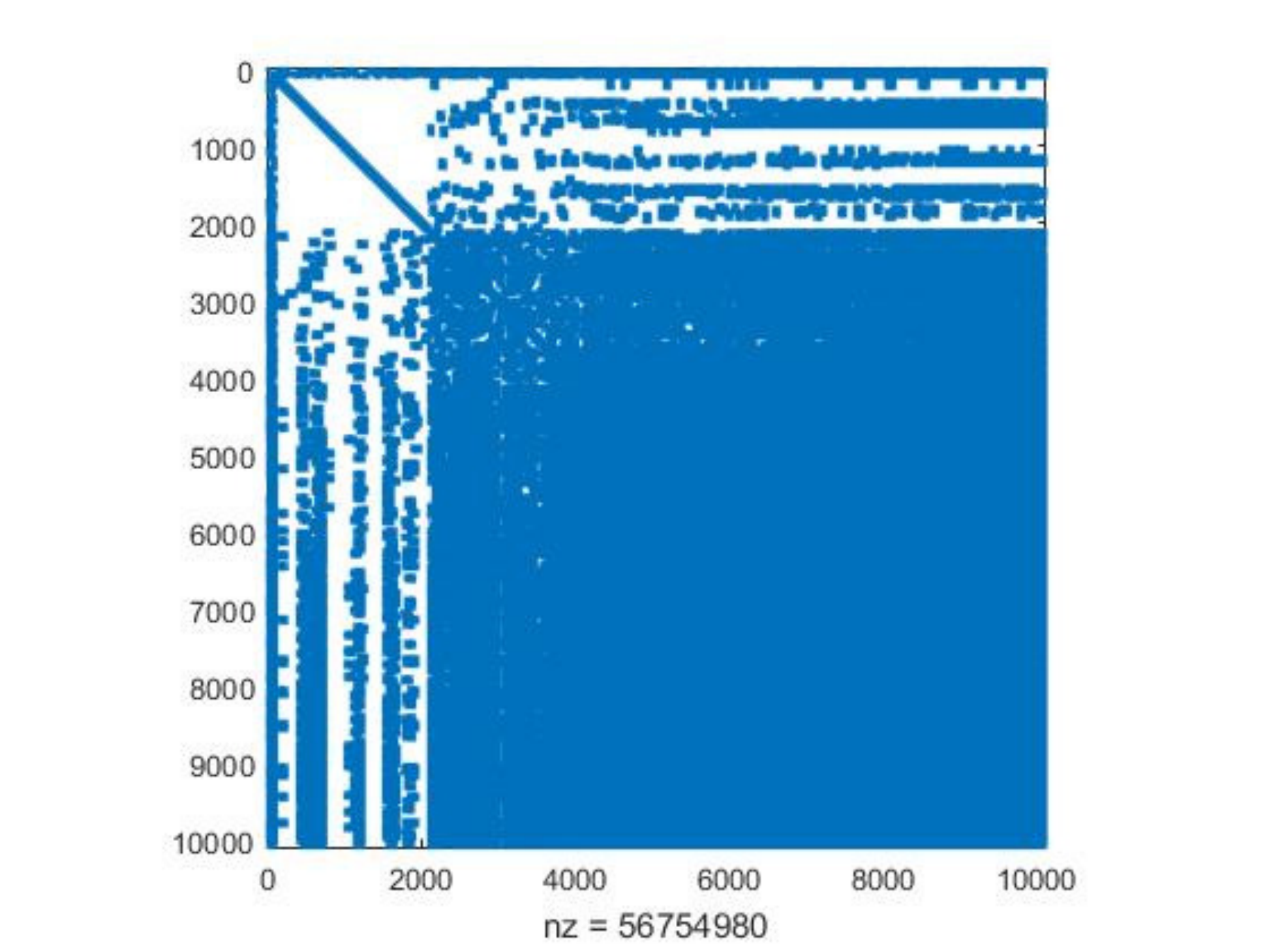}
\end{center}
\end{minipage}
\hfill
\begin{minipage}{.5\textwidth}
\begin{center}
\includegraphics[width=200pt]{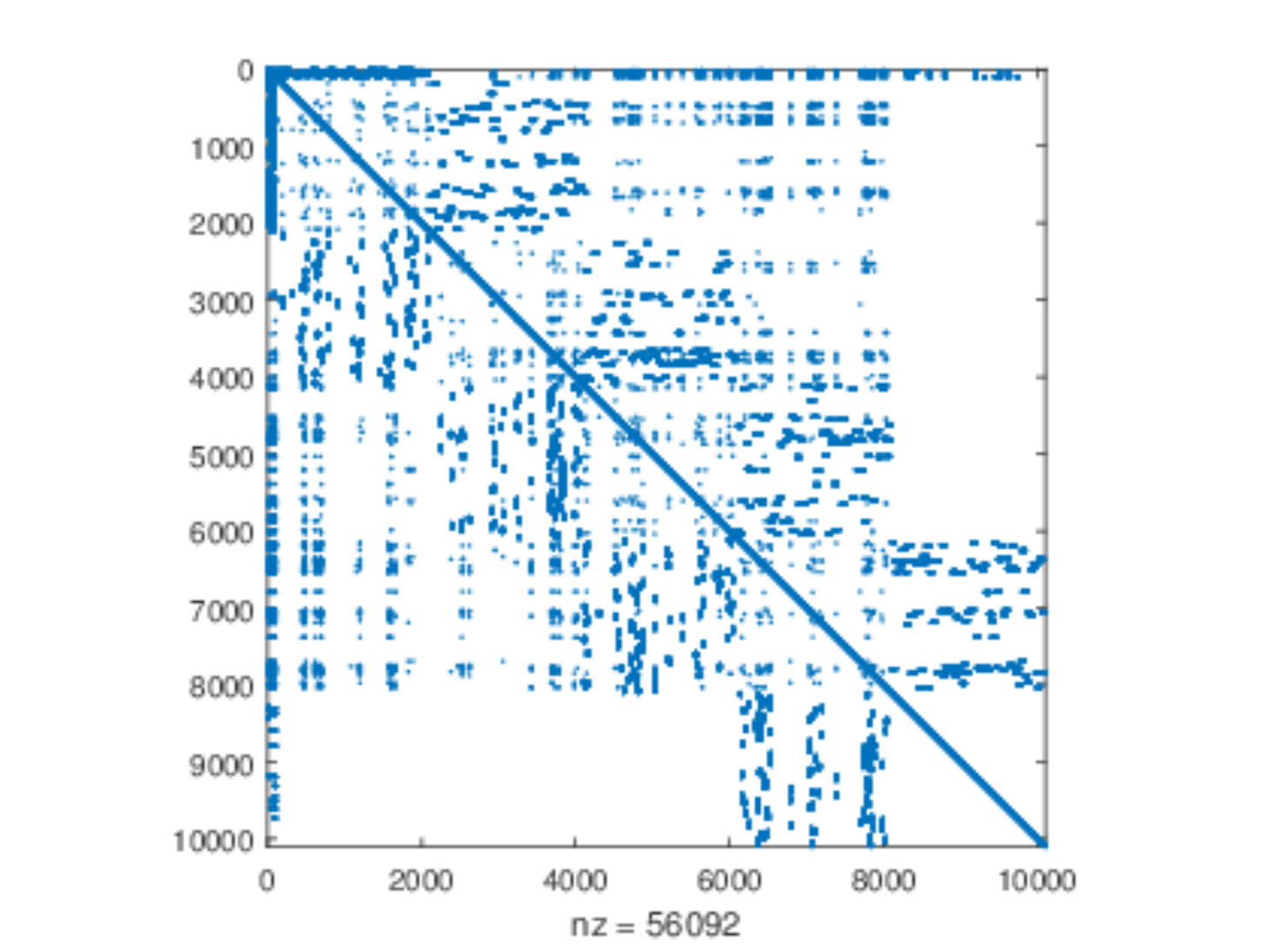}
\end{center}
\end{minipage}
\caption{The positions of non-zero elements of $\A$ (left) and $\A^{-1}$ (right) 
for the problem of size $m=10,100$.}\label{fig:AandAinv}
\end{figure}

To obtain an efficient formulation based on SOCP, we investigated key properties 
of the simple SOCP formulation~(\ref{eq:straightSOCP}).
In particular, we focused on the structure of the numerator relationship
matrix $\A$ and its inverse $\A^{-1}$, 
since the SDP formulation uses only $\A^{-1}$.
Figure~\ref{fig:AandAinv} illustrates the positions of the non-zero elements of $\A$ and $\A^{-1}$
for the problem of size $m=10,100$.
The dimensions of $\A$ and $\A^{-1}$ corresponds to the size $m$.
We can see from Figure~\ref{fig:AandAinv} that 
$\A^{-1}$ is much sparser than $\A$.
The numbers of non-zero elements in $\A$ and $\A^{-1}$ are 56,754,980 and 56,092, respectively, hence
their density against the fully-dense matrix ($m^2$ non-zero elements) are
55.6\% and 0.0549\%, respectively.
When we apply the Cholesky factorization to $\A$, the upper-triangular matrix $\U$
inherits the dense property. The number of non-zero elements in $\U$ is 23,171,296,
and its density against the fully-dense upper triangular matrix is 45.4\%. 
We should utilize the property that $\A^{-1}$ is remarkably sparse compared to $\A$ and $\U$.

From this observation, the direction we should pursue is to use $\A^{-1}$ instead of $\A$ itself.
A key step of our approach is
to introduce the new variable $\y = \A \x$. Then, the replacement
$\x$ by $\A^{-1} \y$ in the optimal selection~(\ref{eq:OPSEL}) leads to
an equivalent optimization problem:
\begin{eqnarray*}
\begin{array}{lcl}
\max &:& (\A^{-1}\g)^T \y \\
\mbox{subject to} &:& (\A^{-1}\e)^T \y = 1 \\ 
& & \frac{\y^T \A^{-1} \y}{2} \le \theta \\
& & \l \le \A^{-1}\y \le \u.
\end{array} 
\end{eqnarray*}
The Cholesky factor of $\A^{-1}$ is the transposed matrix of $\U^{-1}$, denoted by $\U^{-T}$:
namely, $\A^{-1} = (\U^{-T})^T \U^{-T}$.
In practical implementation, since $\A^{-1}$ is remarkably sparse,
we apply an appropriate row/column 
permutation like AMD
(approximate minimum degree permutation)~\cite{AMESTOY04}
to $\A^{-1}$ 
in order to reduce the number of fill-in (the non-zero elements 
that newly appear in $\U^{-T}$ during the process of the Cholesky factorization).
Using $\U^{-T}$, we convert this problem into an SOCP problem:
\begin{eqnarray}
\begin{array}{lclcrcc}
\max &:& & & (\A^{-1}\g)^T \y \label{eq:sparseSOCP} \\
\mbox{subject to} &:& 
\left(\begin{array}{r}
1 \\ -1 \\ \u \\ -\l \\
\hline
\sqrt{2 \theta} \\ \0
\end{array} \right)
&-&
\left(\begin{array}{r}
(\A^{-1}\e)^T \\ -(\A^{-1}\e)^T \\ \A^{-1} \\ -\A^{-1} \\ 
\hline
0 \\ \U^{-T}
\end{array} \right) \y & \in &
\Real_+^{2 + 2m} \times \KC^{1+m},
\end{array}
\end{eqnarray}
We call this formulation {\it a sparse SOCP formulation} of the optimal selection 
problem~(\ref{eq:OPSEL}). 

We use the matrix $\A$ to formulate this new SOCP formulation,
but we do not have to use $\A$ to solve the resultant SOCP problem with the interior-point methods.
Furthermore, when we obtain the optimal solution $\y^*$ of~(\ref{eq:sparseSOCP}),
we can also obtain
the optimal solution $\x^*$ of the original problem~(\ref{eq:OPSEL}) 
via the relation $\x^* = \A^{-1} \y^*$ without using the dense matrix $\A$.

In Table~\ref{table:sparseSOCP},  we compare the performance of 
the SDP formulation~(\ref{eq:SDP}), the simple SOCP formulation~(\ref{eq:straightSOCP}),
and the sparse SOCP formulation~(\ref{eq:sparseSOCP}).
In this table, the first column `nnz' is the total number of non-zero elements 
of $\F_0, \F_1, \ldots, \F_m$ for SDP and that of $\f_0$ and $\F$ for SOCP.
The second column is the computation time to convert 
the pedigree like Figure~\ref{fig:tree} and EBVs (estimated breeding values) 
into the SDP or SOCP formulations, and
the third column is the computation time of the solvers.
We applied SDPA with four cores to the SDP formulation and ECOS to the SOCP formulations.
The fourth (last) column is the total computation time.
\begin{table}
\caption{Performance comparison on the 
the SDP formulation, the simple SOCP formulation,
and the sparse SOCP formulation (time in seconds).}
\label{table:sparseSOCP}
\begin{center}
\begin{tabular}{r|rrrr}
\hline
\multicolumn{5}{c}{$m$ (size of pedigree) = 2,045} \\
\hline
 & \multicolumn{1}{|r}{nnz} & time (conversion) & time (solver) & time (total) \\
\hline
SDP formulation (\ref{eq:SDP})& 24300 & 0.52 & 69.55 & 70.21 \\
simple SOCP formulation (\ref{eq:straightSOCP}) & 18201 & 0.10 & 0.04 & 0.28 \\
sparse SOCP formulation (\ref{eq:sparseSOCP})  & 30348 & 0.37 & 0.05 & 0.55 \\
\hline
\multicolumn{5}{c}{$m$ (size of pedigree) = 10,100} \\
\hline
 & \multicolumn{1}{|r}{nnz} & time (conversion) & time (solver) & time (total) \\
\hline
SDP formulation (\ref{eq:SDP}) & 121703 & 26.97 & 39173.03 & 39200.78 \\
simple SOCP formulation (\ref{eq:straightSOCP}) & 23231801 & 15.95 & 5587.53 & 5604.25 \\
sparse SOCP formulation (\ref{eq:sparseSOCP})  & 159570 & 24.14 & 0.68 & 25.60 \\
\hline
\end{tabular}
\end{center}
\end{table}

Table~\ref{table:sparseSOCP} shows that for the large problem, the computation time of the sparse SOCP formulation 
was much shorter than that of the simple SOCP formulation.
The sparse SOCP formulation seemed to have more complex structure 
than the simple SOCP formulation since it repeatedly contained $\A^{-1}$ in the matrix $\F$, but
the total number of non-zeros was reduced from 23,231,801 in the simple SOCP formulation to
159,570 in the sparse SOCP formulation. This led the computation time reduction for the SOCP solver.

To reduce the total time further, we now address the computation time to build the SOCP formulation.
In the case of $ m = 10,100$, the conversion time occupied 94~\% of the total time.
In particular, the construction of the dense matrix $\A$ and its inversion are the principal bottlenecks.

We investigate further the properties of $\A^{-1}$, and employ a compact algorithm to construct
$\A^{-1}$ proposed by Henderson~\cite{HENDERSON76}.
In the compact algorithm, we use the vector of inbreeding coefficients
defined by $\h := diag(\A) - \e$, where  
$diag(\A)$ is a vector composed of the diagonal elements of $\A$.
Quaas~\cite{QUAAS76} devised an efficient method to compute the inbreeding coefficients $\h$ 
without constructing the matrix $\A$ itself,
and Masuda~{\it et al.} utilized this method to implement their YAMS 
package~\cite{masuda2014application}.
The compact algorithm in~\cite{HENDERSON76} 
with the enhancement~\cite{QUAAS76} is summarized in Table~\ref{table:compact}.
We use an indicator function $\delta(p)$ such that $\delta(0) = 1$, and $\delta(p) = 0$ for $p \ne 0$.
We also use the notation $A_{ij}^{-1}$ to denote the $(i,j)$ element of $\A^{-1}$.

\begin{table}
\caption{A compact algorithm to obtain the inverse of the numerator relationship matrix}
\label{table:compact}
\noindent \hspace{2cm} $\A^{-1} \leftarrow \O$. \newline
\noindent \hspace*{2cm}  {\bf for} $ i = 1, 2, \ldots, m$ \newline
\noindent \hspace*{3cm}  $b_i \leftarrow \frac{4}{(1+\delta(p(i))(1-h_{p(i)}) + (1+\delta(q(i))(1-h_{q(i)})}$ \newline
\noindent \hspace*{3cm} {\bf if} $ i \in \PC_0$ {\bf then} \newline
\noindent \hspace*{4cm} add $b_i$ to $A_{ii}^{-1}$ \newline
\noindent \hspace*{3cm} {\bf elseif} $ i \in \PC_1$ {\bf then} \newline
\noindent \hspace*{4cm} add $b_i$ to $A_{ii}^{-1}$ \newline
\noindent \hspace*{4cm} add $-\frac{b_i}{2}$ to $A_{i,p(i)}^{-1}$, and
$A_{p(i),i}^{-1}$ \newline
\noindent \hspace*{4cm} add $\frac{b_i}{4}$ to $A_{p(i),p(i)}^{-1}$ \newline
\noindent \hspace*{3cm} {\bf elseif} $ i \in \PC_2$ {\bf then} \newline
\noindent \hspace*{4cm} add $b_i$ to $A_{ii}^{-1}$ \newline
\noindent \hspace*{4cm} add $-\frac{b_i}{2}$ to $A_{i,p(i)}^{-1}$ ,
$A_{p(i),i}^{-1}$,  $A_{i,q(i)}^{-1}$, and  $A_{q(i),i}^{-1}$ \newline
\noindent \hspace*{4cm} add $\frac{b_i}{4}$ to $A_{p(i),p(i)}^{-1}$ ,
$A_{p(i),q(i)}^{-1}$,  $A_{q(i),p(i)}^{-1}$, and $A_{q(i),q(i)}^{-1}$ \newline
\noindent \hspace*{3cm} {\bf endif} \newline
\noindent \hspace*{2cm}  {\bf endfor} 
\end{table}

The compact formula for Table~\ref{table:compact} can be described as
\begin{eqnarray*}
\A^{-1} &=& \sum_{i \in \PC_0} b_i
\e_i \e_i^T 
+ \sum_{i \in \PC_1} b_i
\left(\e_i - \frac{1}{2} \e_{p(i)} \right) 
\left(\e_i - \frac{1}{2} \e_{p(i)} \right)^T  \\
& & + \sum_{i \in \PC_2} b_i
\left(\e_i - \frac{1}{2} \e_{p(i)} -  \frac{1}{2} \e_{q(i)}\right) 
\left(\e_i - \frac{1}{2} \e_{p(i)} - \frac{1}{2} \e_{q(i)}\right)^T.
\end{eqnarray*}
When we apply this formula to the example in Figure~\ref{fig:tree},
we obtain the inverse of the numerator relationship matrix as follow:
{\small
\begin{eqnarray*}
\A^{-1} = 
\frac{1}{42}
\left(\begin{array}{rrrrrrrrr}
105  &  42 &  -42  & -42 &   21  &   0  & -42 &    0  &   0 \\
    42  &  98 &  -42 &  -42 &  -28 &    0  &   0  &   0  &   0 \\
   -42  & -42 &  105 &   21 &    0  & -42   &  0  &   0 &    0 \\
   -42 &  -42 &   21  & 105  &   0  & -42  &   0   &  0  &   0 \\
    21  & -28  &   0  &   0 &   98 &    0  & -21   &  0 &  -42 \\
     0  &   0 &  -42  & -42  &   0  & 108  &  24  & -48 &    0 \\
   -42  &   0  &   0  &   0  & -21  &  24  & 129  & -48 &  -42 \\
     0  &   0  &   0  &   0  &   0  & -48  & -48 &   96 &    0 \\ 
     0  &   0  &   0  &   0  & -42  &   0  & -42  &   0 &   84
\end{array}\right).
\end{eqnarray*}
} 
We can see in this example that $\A^{-1}$ has more zero-elements than $\A$ of (\ref{eq:A}).
 
It is known that $b_i > 0$ for any $i = 1, 2, \ldots, m$ from properties of the inbreeding coefficients.
Therefore, the constraint on the group coancestry can be transformed into 
another second-order cone constraint:
\begin{eqnarray}
& & \frac{\x^T \A \x}{2} \le \theta 
\quad \Leftrightarrow \quad \frac{\y^T \A^{-1} \y}{2} \le \theta \nonumber \\
& \Leftrightarrow &   \sum_{i \in \PC_0} b_i
y_i^2
+ \sum_{i \in \PC_1} b_i
\left(y_i - \frac{1}{2} y_{p(i)} \right)^2 
+ \sum_{i \in \PC_2} b_i
\left(y_i - \frac{1}{2} y_{p(i)} -  \frac{1}{2} y_{q(i)}\right)^2 
\le 2 \theta \nonumber
\\
& \Leftrightarrow & 
\left(\begin{array}{c}
\sqrt{2 \theta} \\ \B \y 
\end{array}\right) \in
\KC^{1+m}. \label{eq:simpleSOCP}
\end{eqnarray}
Here,  $\B$ is the matrix whose $i$th row vector is defined by
\begin{eqnarray*}
B_{i*} = \left\{\begin{array}{lcl}
\sqrt{b_i} \e_i^T & \mbox{for} & i \in \PC_0 \\
\sqrt{b_i} \left(\e_i - \frac{1}{2} \e_{p(i)} \right)^T & \mbox{for} & i \in \PC_1 \\
\sqrt{b_i} \left(\e_i - \frac{1}{2} \e_{p(i)} - \frac{1}{2} \e_{q(i)} \right)^T & \mbox{for} & i \in \PC_2.
\end{array}\right.
\end{eqnarray*}

We now replace the matrix $\U^{-T}$ in 
the sparse SOCP formulation~(\ref{eq:sparseSOCP})
by $\B$, hence, we derive another SOCP formulation:
\begin{eqnarray}
\begin{array}{lclcrcc}
\max &:& & & (\A^{-1}\g)^T \y \label{eq:compactSOCP} \\
\mbox{subject to} &:& 
\left(\begin{array}{r}
1 \\ -1 \\ \u \\ -\l \\ \hline 
\sqrt{2 \theta} \\ \0
\end{array} \right)
& - &
\left(\begin{array}{r}
(\A^{-1}\e)^T \\ -(\A^{-1}\e)^T \\ \A^{-1} \\ -\A^{-1} \\ 
\hline 
0 \\ \B
\end{array} \right) \y &\in &
\Real_+^{2 + 2m} \times \KC^{1+m},
\end{array}
\end{eqnarray}
We call this formulation {\it a compact SOCP formulation}.

The combination of 
the compact algorithm in Table~\ref{table:compact} and
the second-order cone constraint~(\ref{eq:simpleSOCP}) 
gives us a formulation that does not rely on any dense matrices.
Table~\ref{table:compactSOCP} adds the results of 
the compact SOCP formulation to Table~\ref{table:sparseSOCP}.
From Table~\ref{table:compactSOCP}, we observe that
the solver time for the compact SOCP formulation 
was slightly shorter than the sparse SOCP formulation.
A further computation time reduction was obtained in the computation time
to build the SOCP formulations from the pedigree.
In the case $m = 10,100$, the conversion was reduced
from 24.14 seconds to 0.37 seconds.
Furthermore, we saved a lot of memory space. 
For the case $m = 10,100$, the matrix $\A$ required 778 MB of memory,
while in contrast  the memory required for $\B$ is only 2.33 MB.

\begin{table}
\caption{Performance comparison on the 
the SDP formulation and the SOCP formulations
(time in seconds).}
\label{table:compactSOCP}
\begin{center}
\begin{tabular}{r|rrrr}
\hline
\multicolumn{5}{c}{$m$ (size of pedigree) = 2,045} \\
\hline
 & \multicolumn{1}{|r}{nnz} & time (conversion) & time (solver) & time (total) \\
\hline
SDP formulation (\ref{eq:SDP})& 24300 & 0.52 & 69.55 & 70.21 \\
simple SOCP formulation (\ref{eq:straightSOCP}) & 18201 & 0.10 & 0.04 & 0.28 \\
sparse SOCP formulation (\ref{eq:sparseSOCP})  & 30348 & 0.37  & 0.05 & 0.55 \\
compact SOCP formulation (\ref{eq:compactSOCP})  & 30249  & 0.01 & 0.05 & 0.19 \\
\hline
\multicolumn{5}{c}{$m$ (size of pedigree) = 10,100} \\
\hline
 & \multicolumn{1}{|r}{nnz} & time (conversion) & time (solver) & time (total) \\
\hline
SDP formulation (\ref{eq:SDP}) & 121703 & 26.97 & 39173.03 & 39200.78 \\
simple SOCP formulation (\ref{eq:straightSOCP}) & 23231801 & 15.95 & 5587.53 & 5604.25 \\
sparse SOCP formulation (\ref{eq:sparseSOCP})  & 159570 & 24.14 & 0.68 & 25.60 \\
compact SOCP formulation (\ref{eq:compactSOCP})  & 153001 & 0.37 & 0.62 & 1.76 \\
\hline
\end{tabular}
\end{center}
\end{table}

\section{Numerical tests}\label{sec:test}

We conducted a numerical evaluation of the SDP and SOCP formulations on several datasets.
The datasets are for problems of sizes 2045, 5050, 15100, 15222, 50100, 100100 and 300100.
The data with the sizes 2,045 and 15,222 were from Scots pine orchards and loblolly pine orchards,
respectively, and these data are available at the Dryad Digital Repository
\url{http://dx.doi.org/10.5061/dryad.9pn5m}.
The other data were 
generation by simulation of five cycles of breeding in a closed population using the approach of
\cite{mullin2010using, mullin1995stochastic}.

Although we used Matlab for the comparison between the formulations,
we also implemented the compact SOCP formulation with C++.
There were two reasons to implement it outside a Matlab environment.
The first is that we can directly know the structure of  non-zero elements that appear in the 
matrix $\B$. Therefore, a specified data structure can accelerate the computation time to 
arrange the input data for building $\F$.
Secondly, 
we expect a software package that does not depend on commercial software
would extend the opportunity for the field application in tree breeding.
Due to the latter reason, we also excluded commercial SOCP solvers from the numerical tests.

One of the advantages of the compact SOCP formulation is that we do not need numerical routines
for the Cholesky factorization.
If we implement 
the simple or sparse SOCP formulation with C++, we have to embed certain numerical routines
for the Cholesky factorization. In addition, the sparse Cholesky factorization requires 
a preprocessing by AMD~\cite{AMESTOY04} to derive its best performance.
In contrast, the compact SOCP formulation obtains the matrix directly $\F$ as discussed in 
Section~\ref{sec:SOCP}.

The computation environment for the small problems ($m \le 10,100$) was
Matlab R2015a on a Windows 8.1 PC with Xeon E3-1231 (3.40 GHz) and 8 GB memory space.
For the large problems ($m \ge 15,100$), we used
Matlab R2014b on a Debian Linux server with Opteron 4386 (3.10 GHz) and 128 GB memory space,
since 8 GB memory space was not enough for the SDP formulation.

Table~\ref{table:results} shows the numerical results of the SDP and SOCP formulations.
We observe from this table that the SDP formulation demanded rather long computation time, particularly
 for the larger problems.
For the problem $m=15,222$, the compact SOCP formulation with C++ reduced the $22,566$ seconds 
of the SDP formulation to only $2.21$ seconds.

\begin{table}
\caption{Numerical results on SDP and SOCP formulations
(time in seconds).}
\label{table:results}
\begin{center}
{\small
\begin{tabular}{r|rrrr}
\hline
\multicolumn{5}{c}{$m$ (size of pedigree) = 2,045} \\
\hline
 & \multicolumn{1}{|r}{nnz} & time (conversion) & time (solver) & time (total) \\
\hline
SDP formulation (\ref{eq:SDP})& 24300 & 0.52 & 69.55 & 70.21  \\
simple SOCP formulation (\ref{eq:straightSOCP}) & 18201 & 0.10 & 0.04 & 0.28 \\
sparse SOCP formulation (\ref{eq:sparseSOCP})  & 30348 & 0.37 & 0.05 & 0.55\\
compact SOCP formulation (\ref{eq:compactSOCP})  & 30249 & 0.01 & 0.05 & 0.20\\
compact SOCP formulation with C++ & 28246 & 0.01 & 0.06 & 0.09\\
\hline
\multicolumn{5}{c}{$m$ (size of pedigree) = 5,050} \\
\hline
 & \multicolumn{1}{|r}{nnz} & time (conversion) & time (solver) & time (total) \\
\hline
SDP formulation (\ref{eq:SDP})& 60853 & 4.63 & 887.32 & 892.30 \\
simple SOCP formulation (\ref{eq:straightSOCP}) & 6812127 & 1.60 & 696.10 & 698.15  \\
sparse SOCP formulation (\ref{eq:sparseSOCP})  & 78405 & 3.67 & 0.19 & 4.21  \\
compact SOCP formulation (\ref{eq:compactSOCP})  & 76533 & 0.04 & 0.19 & 0.58 \\
compact SOCP formulation with C++ & 76533 & 0.01 & 0.21 & 0.28\\
\hline
\multicolumn{5}{c}{$m$ (size of pedigree) = 15,100} \\
\hline
 & \multicolumn{1}{|r}{nnz} & time (conversion) & time (solver) & time (total) \\
\hline
SDP formulation (\ref{eq:SDP}) & 181703 & 157.41 & 21836.63 & 21994.87 \\
simple SOCP formulation (\ref{eq:straightSOCP}) & 54063065  & 26.17  & 38733.01 & 38760.00 \\
sparse SOCP formulation (\ref{eq:sparseSOCP})  & 234760 & 145.53 & 2.13 & 148.49 \\
compact SOCP formulation (\ref{eq:compactSOCP})  & 227989 &  0.04 & 2.06 & 2.92 \\
compact SOCP formulation with C++  & 227989 & 0.02 & 1.95 & 1.99 \\
\hline
\multicolumn{5}{c}{$m$ (size of pedigree) = 15,222} \\
\hline
 & \multicolumn{1}{|r}{nnz} & time (conversion) & time (solver) & time (total) \\
\hline
SDP formulation (\ref{eq:SDP}) & 181947 & 161.99 & 22403.30 & 22566.11 \\
simple SOCP formulation (\ref{eq:straightSOCP}) & 7889551 & 17.96 & 618.18 & 636.95  \\
sparse SOCP formulation (\ref{eq:sparseSOCP})  & 227758 & 150.07 & 2.13 & 153.01 \\
compact SOCP formulation (\ref{eq:compactSOCP})  & 227203 & 0.04 & 2.20 & 3.05 \\
compact SOCP formulation with C++  &  227203 & 0.02 & 2.16 & 2.21\\
\hline
\multicolumn{5}{c}{$m$ (size of pedigree) = 50,100} \\
\hline
 & \multicolumn{1}{|r}{nnz} & time (conversion) & time (solver) & time (total) \\
\hline
SDP formulation (\ref{eq:SDP}) &  &  & & OOM${}^*$ \\
simple SOCP formulation (\ref{eq:straightSOCP}) &  &  & & OOM \\
sparse SOCP formulation (\ref{eq:sparseSOCP})  & 759294 & 4989.55 & 7.58 & 4999.90\\
compact SOCP formulation (\ref{eq:compactSOCP})  & 753023 & 0.15 & 7.71 & 10.63\\
compact SOCP formulation with C++ & 753023 & 0.08 & 7.48 & 7.69\\
\hline
\multicolumn{5}{c}{$m$ (size of pedigrees) = 100,100} \\
\hline
 & \multicolumn{1}{|r}{nnz} & time (conversion) & time (solver) & time (total) \\
\hline
SDP formulation (\ref{eq:SDP}) &  &  &  & OOM\\
simple SOCP formulation (\ref{eq:straightSOCP}) &  &  & & OOM \\
sparse SOCP formulation (\ref{eq:sparseSOCP})  &  &  & & $>$ 24 hours${}^{**}$ \\
compact SOCP formulation (\ref{eq:compactSOCP})  & 1502983 & 0.35 & 18.44 & 24.76\\
compact SOCP formulation with C++ & 1502983 & 0.16 & 17.57 & 17.92 \\
\hline
\multicolumn{5}{c}{$m$ (size of pedigrees) = 300,100} \\
\hline
 & \multicolumn{1}{|r}{nnz} & time (conversion) & time (solver) & time (total) \\
\hline
SDP formulation (\ref{eq:SDP}) &  &  & & OOM \\
simple SOCP formulation (\ref{eq:straightSOCP}) &  &  & & OOM \\
sparse SOCP formulation (\ref{eq:sparseSOCP})  &  &  & & OOM\\
compact SOCP formulation (\ref{eq:compactSOCP}) & 4503065 & 1.10 & 82.41 & 99.65 \\
compact SOCP formulation with C++ & 4503065 & 0.51 & 78.54 & 79.62 \\
\hline
\multicolumn{5}{l}{
{\footnotesize
* OOM = ``out of memory"
} 
} 
\\
\multicolumn{5}{l}{
{\footnotesize
* $>$ 24 hours =
``the computation failed to complete with in in 24 hours."
} 
} 
\end{tabular}
} 
\end{center}
\end{table}

Another significant advantage of the compact formulation is memory consumption.
The SDP formulation consumed 31 GB memory space to solve the problem $m = 15,222$,
and it failed to handle $m \ge 50,100$, despite having  128 GB of memory available.
If we use a rough estimation, the memory required for the largest problem $m = 300,100$ would be 12,000 GB.
The simple SOCP formulation also suffered from a heavy memory requirement to store the dense matrix $\U$.
The sparse SOCP formulation did not require the dense matrix $\A$ in the resultant SOCP problem,
but it utilized $\A$ and computed $\A^{-1}$, hence it required the long conversion time in the same way as the SDP formulation.
In contrast, the compact SOCP formulation consumed less than 766 MB memory space to solve even the largest problem $m = 300,100$.
This memory reduction was mainly a result of employing Henderson's algorithm.

From the numerical results, we also observe that the compact SOCP formulation with C++ is faster than that with Matlab.
The discrepancy between Matlab (99.65 seconds) and C++ (82.41 seconds) in the largest problem 
was due to a specified data structure written in C++.
In particular, the structure was effective when we arranged the pedigree 
before building $\F$. 



\section{Conclusions and future directions}\label{sec:conclusion}

We examined the SOCP formulations for the optimal selection problem arising from 
tree breeding.  We employed the transformation $\x = \A^{-1} \y$ based
on the sparsity of $\A^{-1}$ and
the efficient method to build $\A^{-1}$ by
the Henderson's algorithm with the Quaas enhancement. 
The compact SOCP formulation thus did not involve any dense matrix or the Cholesky factorization.
The numerical results demonstrated that the compact SOCP formulation obtained the optimal solution
significantly faster than the existing SDP formulation.

The SOCP formulations proposed in this paper may look rather simple
for researchers in the field of mathematical optimization.
However,  the computation time reduction in the optimal selection problem
will help improve operational application in tree breeding.
We expect that this paper will be one of bridges to introduce efficient approaches
cultivated in mathematical optimization to tree breeders.

In this paper, we discussed an {\it unequal} deployment of parental genotypes to seed orchards, where
the contributions of selected members are not required to be equal. 
To deal with  {\it equal} deployment, as might be appropriate for selection of a fixed-size breeding population,
we will need a method 
to solve a mixed integer SOCP problem; that is an SOCP problem
 in which some variables are constrained to be integers.
The structure of the SOCP formulation developed in this paper will be a basis 
for an efficient method to consider the mixed integer SOCP problem in the optimal selection problems.
 
\section*{Acknowledgments}

We are grateful to Dr. Yutaka Masuda of 
Obihiro Univeristy of Agriculture and Veterinary Medicine
for providing access to the source code of the YAMS package.
Our work was partially supported by funding from
JSPS KAKENHI (Grant-in-Aid for Scientific Research (C), 15K00032)
and
 F\"{o}reningen Skogstr\"{a}dsf\"{o}r\"{a}dling (The Swedish Tree Breeding Foundation).


\end{document}